\documentclass[11pt]{article} 
\usepackage{amssymb}
\usepackage{amsmath}
\usepackage{algorithmic}
\usepackage{algorithm}

\newenvironment{proof}{\noindent {\bf Proof }}
{\hfill $\bullet$ \vspace{0.25cm}}

\def\one{{\bf 1}\hskip-.5mm}
\def\reff#1{(\ref{#1})}
\def\E{{\mathbb E}}
\def\P{{\mathbb P}}
\def\R{{\mathbb R}}
\def\Z{{\mathbb Z}}

\def\N{{\mathbb N}}
\def\bla{{\bar{\lambda}}}

\newtheorem{theo}{Theorem}
\newtheorem{prop}{\indent Proposition}

\newtheorem{rem}{\indent Remark}
\newtheorem{lem}{\indent Lemma}
\newtheorem{defin}{\indent Definition}
\newtheorem{cor}{\indent Corollary}

\newtheorem{ex}{\indent Example}

\title{Perfect simulation and finitary coding for multicolor systems with interactions of infinite range}

\author{A.~Galves \and N.~L.~Garcia \and E.~L\"ocherbach}

\date{April 4, 2009}

\begin{document}

\maketitle
\begin{abstract}
  We consider a particle system on $\Z^d$ with finite state space and
  interactions of infinite range. Assuming that the rate of change is
  continuous and decays sufficiently fast, we introduce a perfect
  simulation algorithm for the stationary process. The basic tool we
  use is a representation of the infinite range change rates as a
  mixture of finite range change rates. The perfect simulation scheme
  provides the basis for the construction of a finitary coding from an
  i.i.d. finite-value process to the invariant measure of the multicolor
  system.
\end{abstract}



\section{Introduction}

In this paper we first present a perfect simulation
algorithm for a multicolor system on $\Z^d$ with interactions of
infinite range. This perfect simulation algorithm is the basis of the construction of a
finitary coding from a finite-valued i.i.d. process to the invariant
probability measure of the multicolor system.

By a perfect simulation algorithm we mean a simulation which samples
precisely from the stationary law of the process. More precisely, for
any finite set of sites $F$ and any finite time interval $[0, t]$ we
want to sample the stationary time evolution of the coloring of sites
in $F$ during $[0,t]$. 

By coding we mean
a translation invariant deterministic measurable map from the 
 finite-valued i.i.d. process to the invariant probability
measure of the system. Finitary means that the value of the map
at the origin depends only on a finite subset of the random
variables. This finite subset is a function of the realization of the
family of independent random variables.

The process we consider is an interacting particle
system with finite state space. The elements of this finite state
space are called {\sl colors}.  To each site in $\Z^d$ is assigned a
color. The coloring of the sites changes as time goes by. The rate at which
the color of a fixed site $i$ changes from a color $a$ to a new color
$b$ is a function of the entire configuration and depends on $b$.

We do not assume that the system has a dual, or is attractive, or
monotone in any sense.  Our system is not even spatially
homogeneous. The basic assumptions are the continuity of the infinite
range change rates together with a fast decay of the long range
influence on the change rate. These two properties imply that the change rates can be represented as
a countable mixture of local change rates of increasing range. This
decomposition (see Theorem \ref{theo:decomp}) extends to the case of
interacting particle systems the notion of random Markov chains
appearing explicitly in Kalikow (1990) and Bramson and Kalikow (1993)
and implicitly in Ferrari et al. (2000) and Comets et al. (2002).

The decomposition of the change rate of infinite range as a countable
mixture of finite range change rates suggests the construction of any
cylindrical time evolution of the stationary process by the
concatenation of two basic algorithms. First we construct a backward
black and white sketch of the process. Then in a second forward
algorithm we assign colors to the black and white picture.

The proof that the backward black and white algorithm stops after a
finite number of steps follows ideas presented in Bertein and Galves
(1977) to study dual processes. Using these ideas we prove the
existence of our process in a self-contained way. The same ideas
appear again in the construction of the finitary coding. 
This type of construction is similar in spirit to  procedures
adopted in Ferrari (1990), Ferrari et al. (2002), Garcia and Mari\'c
(2006) and Van den Berg and Steif (1999). However all these papers
only consider particular models, satisfying restrictive assumptions
which are not assumed in the present paper.

Our Theorem \ref{theo:5} shows the existence of a finitary coding from
an i.i.d. finite-valued process to the invariant probability measure
of the multicolor system. This can be seen as an extension to the
infinite range processes of Theorem 3.4 of Van den Berg and Steif
(1999). 

H\"aggstr\"om and Steif (2000) constructs a finitary coding for Markov
fields. This result follows, under slightly stronger assumptions, as a
corollary of our Theorem \ref{theo:5} which holds also for non Markov
infinite range fields. These authors conclude the above mentioned
paper by observing that the extension of their results to
``infinite-range Gibbs measures appears to be a more difficult
matter''. Our Theorem \ref{theo:5} is an attempt in this direction.

This paper is organized as follows. In Section 2 we present the model
and state a preliminary result, Theorem \ref{theo:decomp}, which 
gives the representation of the change rate as a countable
mixture of local change rates.  In Section 3, we present the perfect
simulation algorithm and Theorem \ref{theo:nstop} which ensures that
the algorithm stops after a finite number of steps. Theorem
\ref{theo:nstop} also guarantees the exponential ergodicity of the
process. The definitions and results concerning the finitary coding
are presented in Section 4.  The proofs of the theorems are presented
in Sections 5 to 11.

\section{Definitions, notation and basic results}

In what follows, $A$ will be a finite set of colors, the initial
lowercase letters $a$, $b$, $c, \ldots$ will denote elements of $A.$
We will call configuration any element of $A^{\Z^d} .$
Configurations will be denoted by letters  
$\eta, \zeta, \xi , ...$ A point $ i \in \Z^d $ will be 
called site. 
As usual, for any $i \in \Z^d$, $\eta(i)$ will denote the
value of the configuration $\eta$ at site $i$. By extension, for any
subset $V \subset \Z^d$, $\eta(V)\in A^V$ will denote the restriction
of the configuration $\eta$ to the set of positions in $V.$ For any
$\eta ,$ $i$ and $a,$ we shall denote $\eta^{i,a}$ the modified
configuration
$$\eta^{i,a}(j) = \eta(j) \mbox{, for all $j \neq i,$ and $\eta^{i,a}(i)
  = a.$}$$ For any $i \in \Z^d,$ $\eta \in A^{\Z^d}$ and $a \in A,$ 
$a \neq \eta(i),$ and we
denote by $c_i (a, \eta) $ a positive real number. We suppose that 
there exists a constant $\Gamma_i < + \infty $ such that
\begin{equation}\label{eq:boundedrate}
 c_i (a, \eta) \le \Gamma_i ,
\end{equation}
for every $a $ and $ \eta$ such that $a \neq \eta(i).$   

A multicolor system with interactions of infinite range is a Markov
process on $A^{\Z^d} $ whose generator is defined on cylinder functions by
\begin{equation}
  \label{eq:generator}
  L \,f(\eta) \,=\, \sum_{i \in \Z^d} \sum_{a \in A, a \neq \eta(i)} 
c_i ( a, \eta) [f(\eta^{i,a}) - f(\eta)]  \, .
\end{equation}
Intuitively, this form of the generator means that the site $i$ will
be updated to the symbol $a,$ $ a \neq \eta(i),$ at a rate $c_i (a, \eta ) $ whenever the
configuration of the system is $\eta .$ The choice of $c_i (a, \eta) $ for
$ a = \eta(i) $ does not affect the generator (\ref{eq:generator}) and
represents a degree of freedom in our model. In what follows we choose $c_i (\eta(i), \eta) $ 
in such a way that 
\begin{equation}\label{eq:rateprobability}
 c_i (a, \eta) = M_i \, p_i (a| \eta) .
\end{equation}
In the above formula $ M_i < + \infty $ is a suitable constant and for every fixed configuration $\eta
,$ $ p_i (\cdot | \eta) $ is a probability measure on $A.$ Condition
(\ref{eq:boundedrate}) implies that such a choice is always possible, for instance by taking $M_i =
|A| \, \Gamma_i $ and defining
\begin{equation}\label{eq:fixedchoice}
c_i ( \eta (i), \eta) = M_i - \sum_{a: a \neq \eta(i)} c_i (a, \eta) .
\end{equation}
We shall call
$(c_i)_{i \in \Z^d} $ a family of rate functions for this fixed choice 
(\ref{eq:fixedchoice}).
 
Our first aim is to give sufficient conditions on $c_i(a,\eta)$
implying the existence of a perfect simulation algorithm of the process having generator (\ref{eq:generator}). 
To state these conditions,
we need some extra notation. Let $V_{i} (k) = \{j \in \Z^d; 0 \le \|j
- i \| \le k\},$ where $\|j\| = \sum_{u=1}^d |j_u|$ is the usual
$L_1$-norm of
$\Z^d$. We will impose the following continuity condition on the family of rate functions $c  .$\\
 
{\bf Continuity condition.} For any symbol $a,$ we will assume that
\begin{equation}
  \label{eq:continuity}
\sup_{i \in \Z^d}   \sup_{\eta(V_i(k)) = \zeta(V_i(k))} | c_i (a, \eta) -
  c_i (a, \zeta )| \rightarrow 0 \, ,
\end{equation}
as $k \rightarrow \infty .$ 

Define
\begin{equation}
  \label{eq:alpha0}
  \alpha_{i} (-1) \,=\, \sum_{a \in A} \min \left( \inf_{\zeta \in A^{\Z^d}, \zeta(i ) \neq a} c_i(a, \zeta)
, \, M_i - \sup_{ \zeta \in A^{\Z^d}, \zeta(i )= a} \sum_{b \neq a} c_i (b, \zeta) \right)  \, ,
\end{equation}
and for any $k \ge 0 , $
\begin{equation}
  \label{eq:alpha}
  \alpha_i (k) \,=\, \min_{w \in A^{V_i (k) }} \left( (\sum_{a \in A, a \neq w(i)}
  \inf_{\zeta: \zeta(V_i (k)) = w} c_i( a,\zeta)) + M_i - \sup_{\zeta: \zeta(V_i (k)) = w} \sum_{ b \neq w(i)} c_i (b, \zeta)\right)  .
\end{equation} 

In order to clarify the role of $\alpha_i(k)$, we
present an example that is a spatial version of the chain that
regenerates in $1 :$ the evolution of a site depends on the random
ball around that site where the size of this ball is the smallest
radius such that $1$ belongs to the ball (without the center itself).

\begin{ex}
Let $A = \{ 0, 1 \} , $ 
$$ l_i (\eta) = l, \mbox{ if } \max_{j \in V_i(l), j \neq i } \eta (j) = 0 \mbox{ and } \max_{j
\in V_i (l+1), j \neq i} \eta(j) = 1$$
and define 
$$ c_i(1, \eta ) = q_{l_i (\eta)}, \, c_i(0, \eta ) = 1 - c_i(1, \eta ), $$
where $0 < q_k < 1 $ for all $k.$

Note that in this case,
$$ \sup_i \sup_{\eta(V_i (k)) = \zeta(V_i(k))} | c_i(1, \eta) -
c_i(1, \zeta )| = \sup_{l, m \geq k } |q_l - q_m|,$$ and thus the
process is continuous, if and only if $\lim_k q_k$ exists. 

Observe that if $q_k \downarrow q_{\infty},$ as $k \to \infty ,$ then $M_i = 1$ and
\begin{eqnarray*}
\alpha_i (-1)  &=& \min \left( \inf_{\zeta \in A^{\Z^d}, \zeta(i
    ) = 1} c_i(0, \zeta), \, 1 - \sup_{ \zeta \in A^{\Z^d}, \zeta(i )=
    0} c_i (1, \zeta) \right) \\
&&   + \,  \min \left( \inf_{\zeta \in A^{\Z^d}, \zeta(i) = 0} c_i(1,
  \zeta), \, 1 - \sup_{ \zeta \in A^{\Z^d}, \zeta(i )= 1} c_i (0,
  \zeta) \right)  \\
&=&   \inf_{\zeta \in A^{\Z^d}} c_i(0, \zeta) + \inf_{ \zeta \in A^{\Z^d}} c_i (1,
  \zeta) \\
&=& 1 - q_0 + q_{\infty}.
\end{eqnarray*}
Also,
\begin{eqnarray*}
\alpha_i (0)  &=& \min \Big( \inf_{\zeta \in
    A^{\Z^d}, \zeta(i) = 0} c_i(1, \zeta) + 1 -  \sup_{\zeta \in
    A^{\Z^d}, \zeta(i) = 0} c_i(1, \zeta), \\
&& \quad \quad \inf_{\zeta \in
    A^{\Z^d}, \zeta(i) = 1} c_i(0, \zeta) + 1 -  \sup_{\zeta \in
    A^{\Z^d}, \zeta(i) = 1} c_i(0, \zeta) \Big) \\
&=& \min (q_{\infty} + 1 - q_0, 1 - q_0 + q_{\infty} )\\
&=& \alpha_i (-1) .
\end{eqnarray*}

Finally, 
\begin{eqnarray*}
\alpha_i (k)  &=& \min_{w \in A^{V_i(k)}} \left( \inf_{\zeta:
    \zeta(V_i (k)) = w} c_i( 1 - \zeta(i),\zeta) + 1 - \sup_{\zeta:
    \zeta(V_i (k)) = w}  c_i (\zeta(i), \zeta) \right) \\
&=& \min_{w \in A^{V_i(k)}} \left( \inf_{\zeta:
    \zeta(V_i (k)) = w}  c_i( 0,\zeta) + \inf_{\zeta:
    \zeta(V_i (k)) = w}  c_i( 1,\zeta) \right) \\
&=& \min_{w \in A^{V_i(k)}}  \left( \inf_{\zeta:
    \zeta(V_i (k)) = w} q_{l_i (\zeta)} +  \inf_{\zeta:
    \zeta(V_i (k)) = w} (1-q_{l_i (\zeta)}) \right)\\
&=& (1 - q_k) + q_{\infty} .
\end{eqnarray*}
\end{ex}

Let us introduce some more notation. Note that for each site $i,$
\begin{equation}\label{eq:m}
M_i = \lim_{k \to \infty }
\alpha_i (k)  . 
\end{equation}
Hence to each
site $i$ we can associate a probability distribution $ \lambda_i$ by
\begin{equation}
  \label{eq:lambda0}
  \lambda_{i} (-1)\,=\, \frac{\alpha_{i} (-1)}{M_i} ,
\end{equation}
and for $k \ge 0 $
\begin{equation}
  \label{eq:lambdak}
  \lambda_i (k)\,=\, \frac{\alpha_i (k) \,-\, \alpha_{i} (k-1)}{M_i } .
\end{equation}

We will see that for each $i$ the family of rate functions $c_i(.,.)$
can be represented as a mixture of local rate functions weighted by
$(\lambda_i (k))_{k \geq -1}.$ More formally, we have the following
theorem.

\begin{theo} \label{theo:decomp} Let $(c_i)_{i \in \Z^d} $ be a family
  of rate functions satisfying the conditions (\ref{eq:boundedrate}), (\ref{eq:rateprobability}),
  the continuity condition
  (\ref{eq:continuity}) and the summability condition
  (\ref{eq:condition1}).  Then for any site $i$ there exists a family
  of conditional probabilities $p_i^{[k]}$ on $A$ depending on the
  local configurations $\eta(V_i (k))$ such that
\begin{equation}
\label{cmmc}
c_i (a, \eta) \,=\, M_i \, p_i (a|\eta) , \mbox{ where } p_i(a|\eta) =   \sum_{k \geq -1 }  \lambda_i (k) p_i^{[k]} (a | \eta(V_i (k))).
\end{equation}
As a consequence, the infinitesimal generator $L$ given by
(\ref{eq:generator}) can be rewritten as
\begin{equation}
  \label{eq:generator2}
  L \,f(\eta) \,=\, \sum_{i \in \Z^d} \sum_{a \in A} \sum_{k \ge -1} 
M_i \lambda_i (k) p_i^{[k]}(a| \eta(V_i (k))) [f(\eta^{i,a}) - f(\eta)]\, .
\end{equation}
\end{theo}

\begin{rem}
  Note that for $k = -1,$ $V_i (k) = \emptyset $ and hence $
  p_i^{[-1]} (a | \eta(V_i (k))) = p_i^{[-1]} (a)$ does not depend on
  the configuration.  Therefore, $\lambda_{i} (-1)$ represents the
  spontaneous self-coloring rate of site $i$ in the process.
\end{rem}

The representation given by (\ref{eq:generator2}) provides a clearer
description of the time evolution of the process. We start with an
initial configuration $\eta$ at time zero. This configuration is
updated in a c\`adl\`ag way as follows. For each site $ i \in \Z^d ,$
consider a rate $M_i$ Poisson point process $N^i .$ The Poisson
processes corresponding to distinct sites are all independent. If at
time t, the Poisson clock associated to site $i$ rings, we choose a
range $k$ with probability $\lambda_i (k)$ independently of everything
else. And then, update the value of the configuration at this site by
choosing a symbol $a$ with probability $p_i^{[k]} (a |
\xi^{\eta}_t(V_i(k)))$.

\section{Perfect simulation of the stationary process}
The decomposition (\ref{cmmc}) provided by Theorem \ref{theo:decomp}
suggests an algorithm of perfect simulation for the
multicolor long range interacting system. This is the main
result of this article.  The goal is to sample under
equilibrium the time evolution of any
finite set of sites $F$ during any fixed finite time interval.

We first introduce a simulation procedure to sample 
the time evolution of any
finite set of sites $F$ during any fixed finite time interval $[0, t],$ when starting from
a fixed initial configuration $\eta .$ This simulation
procedure has two stages. First, we draw a backward black and white
sketch in order to determine the set of sites and the succession of
choices affecting the configuration of the set of sites $F$ at time $t.$ 
Then, in the second stage, a forward coloring procedure
assigns colors to every site involved in the black and white
sketch. This will be formally described in Algorithms 1 and 2 below. 

Let us describe the mathematical ideas behind this algorithm. Our goal
is to simulate the configuration of the fixed set of sites $F$ during the
time interval
$[0, t]$ when the process starts from an initial configuration $\eta .$
We climb up the rate $M_j$ Poisson processes $N^j, j \in F ,$ until we find the
last occurrence time before time $t$ where the Poisson clock
rang. Note that the probability that the clock of site $i$ rings first
among all these clocks is given by
  $$ \frac{M_i}{\sum_{j \in  F} M_j } .$$
  Then we have to inspect the configuration at the sites belonging to
  the finite set $V_i (k) $ which is chosen at that time. $V_i (k) $
  is chosen with probability $\lambda_i (k), k \geq -1.$ If $k = -1 $ is chosen, this
  means that the value of $\xi (i) $ at that time is chosen according
  to $p_i^{[-1]} ,$ independently of the other sites, and thus site
  $i$ can be removed from the set $F  .$
 
  Otherwise, if $ k \geq 0,$ we have to include all the sites in $V_i
  (k) $ to the set of sites $F$ and to continue the algorithm. The
  reverse-time checking continues for each point reached previously
  until we find an occurrence time before time $0.$ In this
  case the algorithm stops.
  
  In the second stage, the algorithm assigns colors to all the sites
  that have been involved in the first stage. To begin with, 
  all sites that have not yet chosen a range $-1$, will be colored according
  to the initial configuration $\eta $ at time $0.$   
  And then
  successively, going forwards in time, we assign colors to the
  remaining sites according to $ p^{[ k ]}_i ( \cdot | V_i (k)),$
  where all sites in $V_i (k) $ have already been colored in a
  previous step of the algorithm. Finally we finish with the colors of
  the set of sites $F$ at time $t .$

The finite time simulation Algorithm 1 and 2 uses the following variables.

\begin{itemize}
\item $N$ is an auxiliary variables taking values in the set 
of non-negative integers  
$ \{ 0, 1,2, \ldots \} $
\item  $N_{STOP}$ is a counter taking values in the set of non-negative integers
$ \{ 0, 1, 2, \ldots \} $
\item $T_{STOP} $ is an element of $(0, + \infty )$  
\item
$I $ is variable taking values in $\Z^d$ 
\item
$K$ is a variable taking values in $\{ -1, 0, 1, \ldots \}$
\item
$T$ is an element of $( 0, + \infty ) $
\item
$B = (B_1, B_2, B_3)$ where 
\begin{itemize}
\item
$B_1 $ is an array of
elements of $\Z^d $
\item
$B_2$ is an array of
elements of $\{ -1, 0 , 1 , \ldots \} $
\item
$B_3$ is an array of elements of $ (0 , + \infty ) $
\end{itemize}
\item
$C$ is variable taking values in the set of finite subsets of $\Z^d$ 
\item $ W  $ is an auxiliary variable taking values in $ A$
\item $V$ is an array of elements of $A $
\item $\zeta $ is a function from $\Z^d $ to $A \cup \{ \Delta \} ,$
  where $\Delta $ is some extra symbol that does not belong to $A$

\end{itemize}

\begin{algorithm}[h]
\caption{Backward black and white sketch without deaths} 
\begin{algorithmic}[1]
\STATE {\it Input:} $F$; {\it Output:} $N_{STOP}$, $B$, $ C$, $T_{STOP} $
\STATE  $N \leftarrow 0,$ $N_{STOP} \leftarrow 0 ,$ $ B \leftarrow \emptyset ,$ $ C \leftarrow F, $ $T_{STOP} \leftarrow 0 $ 
\WHILE {$ T_{STOP} < t \mbox{ and } C \neq \emptyset $} 
\STATE Choose a time $T  \in (0, +\infty) $ randomly according to the
exponential distribution with parameter 
$ \sum_{j \in C} M_j  .$ Update 
$$ T_{STOP} \leftarrow  T_{STOP} + T .$$
\STATE $N \leftarrow N+1 .$ 
\STATE Choose a site $I \in C$ randomly according to the distribution 
$$P ( I = i)  = \frac{M_i }{\sum_{j \in C } M_j}$$
\STATE Choose $ K \in \{ -1, 0, 1, \ldots \}$ randomly according to
the distribution
$$P( K  = k) = \lambda_I ( k)$$
\STATE $ C \leftarrow C \cup V_I (K)$ 
\STATE $ B (N) \leftarrow  (I, K, T_{STOP} )$
\ENDWHILE
\STATE $N_{STOP} \leftarrow N $ 
\end{algorithmic}
\end{algorithm}

\begin{algorithm}[h]
\caption{ Forward coloring procedure}
\begin{algorithmic}[1]
\STATE {\it Input:} $N_{STOP}$, $B$, $C$, $\eta (C)$; {\it Output:} $ V$ 
\STATE $ N \leftarrow N_{STOP }   $
\STATE $\zeta(j) \leftarrow \eta (j) $ for all $j \in C;$ $\zeta(j) \leftarrow \Delta $ for all $ j \in  \Z^d \setminus C $
\WHILE {$N \ge 1$} 
\STATE $ (I,K,T) \leftarrow B(N)  .$ 
\IF {$K= -1 $} \STATE Choose $W $ randomly in $A$
according to the probability distribution
$$ P( W = v)  = p_I^{[-1]} (v )$$
\ELSE \STATE{Choose $W $ randomly in $A$
according to the probability distribution
$$ P( W  = v)  = p_I^{[K]} ( v | \zeta ( V_I (K))$$}
\ENDIF
\STATE $ \zeta (I) \leftarrow W $
\STATE $ V (N) \leftarrow W $
\STATE { $ N \leftarrow N-1 $}
\ENDWHILE 
\end{algorithmic}
\end{algorithm}

Using output $V$ of Algorithm 2 
and output $B$ of Algorithm 1 we can
construct the time evolution $(\xi_s (F) , 0 \le s \le t) $ of
the process. This is done as follows.

Denote $I(N), T(N) $ the first and the third coordinate of
the array $B(N)$ respectively. Introduce
the following random times for any $1 \le n \le
N_{STOP},$ 
$$S_n = t - T ( N_{STOP} - n +1 ) .$$ 

\begin{itemize}
\item
For $0 \le s < S_1 $ define $\xi_s (F) = \zeta(F ) .$
\item
For $1 \le n \le N_{STOP},$ for $S_n \le s < S_{n+1}\wedge t ,$ we put 
\begin{itemize}
\item
for all $ i \in F $ such that $i \neq I( N_{STOP} - n + 1 ), $ $\xi_s (i) = \xi_{S_n} (i) ; $
\item
for $ i =   I( N_{STOP} - n + 1 ) ,$ $ \xi_s (i) = V(n) .$ 
\end{itemize}
\end{itemize}

We summarize the above discussion in the following proposition.

\begin{prop} \label{theo:1} Let $(c_i)_{i \in \Z^d}$ be a family of
  continuous rate functions satisfying the conditions of Theorem \ref{theo:decomp}. If
  \begin{equation}
    \label{eq:condition1}
   \sup_{i \in \Z^d} \sum_{k \ge 0} |V_i (k)|    \lambda_i (k)    \,< \, +\infty\,,
  \end{equation}
  then Algorithm 1 stops
  almost surely after a finite number of steps, i.e.
  $$ P( N_{STOP} < + \infty ) = 1 .$$
  Moreover, for any initial configuration $\eta$, there exists a unique
  Markov process $(\xi^{\eta}_t)_{t \ge 0}$ such that $\xi^\eta_0 = \eta$ and with
  infinitesimal generator 
  \begin{equation}
   \label{eq:generator1}
   L \,f(\eta) \,=\, \sum_{i \in \Z^d} \sum_{a \in A} c_i(a , \eta) [f(\eta^{i,a}) - f(\eta)] \, .
  \end{equation}
  The cylindrical time
  evolution $
  (\xi_s (F) , 0 \le s \le t ) $ simulated in Algorithms 1 and
  2 is a sample from this process $\xi^{\eta} .$
  \end{prop}

We now turn to the main object of this paper, the perfect simulation of
the multicolor long range interacting system under equilibrium. The goal is to sample under
equilibrium the time evolution of any
finite set of sites $F$ during any fixed finite time interval.

We first introduce a simulation procedure to sample from equilibrium
the cylindrical configuration at a fixed time. As before, this simulation
procedure has two stages : First, we draw a backward black and white
sketch in order to determine the set of sites and the succession of
choices affecting the configuration of the set of sites at
equilibrium.  Then, in the second stage, a forward coloring procedure
assigns colors to every site involved in the black and white
sketch. This will be formally described in Algorithms 3 and 4 below. 

The following variables will be used. 

\begin{itemize}
\item $N$ is an auxiliary variables taking values in the set 
of non-negative integers  
$ \{ 0, 1,2, \ldots \} $
\item  $N_{STOP}$ is a counter taking values in the set of non-negative integers
$ \{ 0, 1, 2, \ldots \} $
\item
$I $ is variable taking values in $\Z^d$ 
\item
$K$ is a variable taking values in $\{ -1, 0, 1, \ldots \}$
\item
$B $ is an array of
elements of $\Z^d \times \{ -1, 0 , 1 , \ldots \} $
\item
$C$ is variable taking values in the set of finite subsets of $\Z^d$ 
\item $ W  $ is an auxiliary variable taking values in $ A$
\item $\eta $ is a function from $\Z^d $ to $A \cup \{ \Delta \} ,$
  where $\Delta $ is some extra symbol that does not belong to $A$

\end{itemize}

{\begin{algorithm}[h]
\caption{Backward black and white sketch} 
\label{algo1}
\begin{algorithmic}[1]
\STATE {\it Input:} $F$; {\it Output:} $N_{STOP}$, $B$ 
\STATE  $N \leftarrow 0,$ $N_{STOP} \leftarrow 0 ,$ $ B \leftarrow \emptyset ,$ $ C \leftarrow F, $  
\WHILE {$C \neq \emptyset  $} 
\STATE $N \leftarrow N+1 .$ 
\STATE Choose a site $I \in C$ randomly according to the distribution 
$$P ( I = i)  = \frac{M_i }{\sum_{j \in C } M_j}$$
\STATE Choose $ K \in \{ -1, 0, 1, \ldots \}$ randomly according to
the distribution
$$P( K  = k) = \lambda_I ( k)$$
\IF {$K = -1,$} \STATE {$C \leftarrow C \setminus \{ I \}$} 
\ELSE \STATE $ C \leftarrow C \cup V_I (K)$ \ENDIF
\STATE $ B (N) \leftarrow  (I, K )$
\ENDWHILE
\STATE $N_{STOP} \leftarrow N $
\end{algorithmic}
\end{algorithm}

\begin{algorithm}[h]
\caption{ Forward coloring procedure}
\begin{algorithmic}[1]
\STATE {\it Input:} $N_{STOP}$, $B$; {\it Output:} $\{(i,\eta(i)),
i \in F\}$ 
\STATE $ N \leftarrow N_{STOP }   $
\STATE $\eta(j) \leftarrow \Delta $ for all $ j \in  \Z^d  $
\WHILE {$N \ge 1$}
\STATE $ (I,K) \leftarrow B(N)  .$ 
\IF {$K= -1 $} \STATE Choose $W $ randomly in $A$
according to the probability distribution
$$ P( W = v)  = p_I^{[-1]} (v )$$
\ELSE \STATE{Choose $W $ randomly in $A$
according to the probability distribution
$$ P( W  = v)  = p_I^{[K]} ( v | \eta ( V_I (K))$$}
\ENDIF
\STATE $ \eta (I) \leftarrow W $
\STATE { $ N \leftarrow N-1 $}
\ENDWHILE 
\end{algorithmic}
\end{algorithm}
}

Let us call $\mu $ the distribution on $A^{\Z^d} $ whose projection on 
$A^F $ is the law of $\eta (F)$ printed at the end of
Algorithm 4. 
The following theorem gives a sufficient condition
ensuring that Algorithm 3 stops after a finite number of steps
and shows that $\mu $ is actually the invariant measure of the process.

\begin{theo}\label{theo:nstop}
Let $(c_i)_{i \in \Z^d} $ be a family
  of rate functions satisfying the conditions of Theorem 
  \ref{theo:decomp}. If  
  \begin{equation}
    \label{eq:condition2}
 \sup_{i \in \Z^d}   \sum_{k \ge 0} \, |V_i (k)|  \lambda_i (k)  \,< \, 1\,,
  \end{equation}
then
$$ P ( N_{STOP} < + \infty ) = 1.$$
The law of the set $\{ (i,\eta (i)) : i \in F \}   $ printed at the
end of Algorithms 3 and 4 is the projection on  $A^F$ of the unique 
invariant probability measure $\mu $ 
of the process. Moreover, the law of the process starting
  from any initial configuration 
  converges weakly to $\mu $ and this convergence takes place
exponentially fast.    
\end{theo}

\begin{rem}
  In the literature, we say that the process is {\it ergodic}, if it
  admits a unique invariant
  measure which is the weak limit of the law of the process starting
  from any initial configuration. If this convergence takes place exponentially fast,
  we say that the process is {\it exponentially ergodic}. Therefore, Theorem
  \ref{theo:nstop} says that the multicolor system is exponentially ergodic.
\end{rem}

Algorithms 3 and 4 show how to sample the invariant probability measure of 
the process. We  now pursuit a more ambitious goal : how to sample
the stationary time evolution of any fixed finite set of sites $F$
during any fixed interval of time $ [0, t ] .$ This is done using 
Algorithms 1 and 2 as well.
 
Algorithm 1 produces a backward black and white sketch without
removing the spontaneously coloring sites.  We start at time $t$ with
the set of sites $F$ and run backward in time until time $0 .$ This
produces as part of its output the set of sites $C$ whose coloring at
time $0$ will affect the coloring of the sites in $F$ during $[0, t].$
We then use the output set $C$ of Algorithm 1 as input set of
positions in Algorithms 3 and 4. Algorithms 3 and 4 will give us as
output the configuration $\eta (C)$ that will be used as input
configuration for Algorithm 2.\\

\begin{theo}\label{theo:4}
Under the conditions of Theorem \ref{theo:1}, Algorithm 3 stops
  almost surely after a finite number of steps
$$ P( N_{STOP} < + \infty ) = 1 .$$
Moreover, under the conditions of Theorem \ref{theo:nstop}, for any
$ t> 0 ,$ the cylindrical time
evolution $
(\xi_s (F) , 0 \le s \le t ) $ simulated in Algorithms 1, 2, 3 and
4 is a sample from the stationary process.
\end{theo}
 
\section{Finitary coding}
The perfect simulation procedure described in Algorithms 1--4 gives
the basis for the construction a finitary coding for the invariant
probability measure of the multicolor system $\xi_t $. By this we mean
the following. Let $(Y(i) , i \in \Z^d)$ be a family of i.i.d. random
variables assuming values on a finite set $S$. Let $ (\xi_0 (i), i \in
\Z^d ) $ be the configuration sampled according to the invariant
probability measure $\mu $ obtained as output of Algorithm 2.

\begin{defin} 
  We say that there exists a {\it finitary coding} from $(Y(i) ,i \in
  \Z^d) $ to $(\xi_0 (i) , i \in \Z^d) $ if there exists a
  deterministic function $f: S^{\Z^d} \rightarrow A^{\Z^d}$
  such that almost surely the following holds:
\begin{itemize}
\item $f$ commutes with the shift operator, that is, $f(T_i(y)) =
  T_i(f(y))$ for any $i \in \Z^d$;
\item $\xi_0 = f( (Y (j) ), j \in \Z^d ) $; and
\item there exists a a finite subset $\bar{F} $ of $\Z^d$ satisfying 
$$  f ( (Y (j) ), j \in \Z^d) =  f ( (Y' (j) ), j \in \Z^d)$$
whenever
$$ Y' (j) = Y (j) \mbox{ for all } j \in \bar{F} .$$
\end{itemize}
\end{defin}

In the first condition of the definition, the notation $T_i$ 
denotes the translation by $i$ steps, in $S^{\Z^d}$,
or in $A^{\Z^d}$. More precisely, for any $i \in \Z^d$, if $y \in
S^{\Z^d}$ then $T_i(y)$ is the element of $S^{\Z^d}$ such that
$T_i(y)(j) = y(j-i)$ with equivalent definition for $\xi \in
A^{\Z^d}$.

\begin{theo}\label{theo:5}
  Under the conditions of Theorem \ref{theo:nstop}
  there exists a finitary coding from an independent and identically
  distributed family of finite-valued variables $ (Y
  (i), i \in \Z^d)$ to $(\xi_0 (i) , i \in \Z^d ) .$
\end{theo}

Theorem 1.1 of H\"aggstr\"om and Steif (2000) follows as a corollary
of Theorem \ref{theo:5} under a slightly stronger condition. In order to
state this corollary, we need to introduce the notion of Markov random field.

\begin{defin} A Markov random field $X$ on $\Z^d$ with values in a
  finite alphabet $A$ has distribution $\mu$ if $\mu$ admits a
  consistent set of conditional probabilities
$$ \mu(X(\Lambda) = \xi(\Lambda)| X(\Z^d \setminus \Lambda) = \xi(\Z^d
\setminus \Lambda))  \,=\,
\mu(X(\Lambda) = \xi (\Lambda)| X(\partial \Lambda) = \xi(\partial
\Lambda))$$
for all finite $\Lambda \subset \Z^d$, $\xi \in A^{\Z^d}$. Here,
$ \partial \Lambda = \{ j \in \Z^d : \inf_{ i \in \Lambda } \| i -j \| = 1 \} .$
Such a set of conditional probabilities is
called the specification of the random field and denoted by ${\cal
  Q}$. 
\end{defin}  

\begin{cor} \label{cor:HS} For any Markov random field $X$ on $\Z^d$ with
  specification ${\cal Q}$ satisfying 
\begin{equation}
\label{eq:HS}
\sum_{a \in A} \min_{\zeta(\partial 0) \in A^{\partial 0}} {\cal Q}(X(0)=a|
X(\partial 0) = \zeta(\partial 0)) >  \frac{2d}{2d +1 },
\end{equation}
there exists an i.i.d. sequence $ (Y (i) , i \in \Z^d ) $ of finite
valued random variables such that there exists a finitary coding from
$ (Y (i) , i \in \Z^d ) $ to the Markov random field.
\end{cor}

\begin{rem}
Just for comparison, Condition (\ref{eq:HS})  is equivalent to
$$\alpha_0(-1)  >  \frac{2d}{2d +1 },$$
while {\it Condition HN} in Theorem 1.1 of H\"aggstr\"om and Steif
(2000)  can be rewritten in our notation as
$$ \alpha_0(-1) > \frac{2d -1 }{2d}.$$ 
This does not seem to be a too high price to pay in
order to be able to treat the general case of long range interactions.
\end{rem}

\section{Proof of Theorem \ref{theo:decomp}}
The countable mixture representation provided by Theorem
\ref{theo:decomp} is the basis of all the other results presented in
this paper. Therefore it is just fair that its proof appears in the
first place.

Recall that $c_i (a, \eta ) = M_i \, p_i (a|\eta ). $ Therefore, it is
sufficient to provide a decomposition for $p_i (a|\eta ).$ Put 
\begin{eqnarray*}
 r_i^{[-1]} (a) &=& \inf_{\zeta} p_i(a |\zeta),\\
   \Delta^{[-1]}_i (a) &= &r_i^{[-1]} (a), \\
r_i^{[0]} ( a| \eta (V_i (0)))& =& \inf_{\zeta : \zeta (V_i (0) ) = 
\eta (V_i (0)) }p_i (a|  \zeta) ,\\
 \Delta^{[0]}_i ( a |  \eta (V_i (0) ))&=
& r_i^{[0]} ( a |\eta(V_i (0) )) - r^{[-1]}_i (a) .
 \end{eqnarray*}
For any $k \geq 1,$ define
$$ r_i^{[k]} ( a| \eta (V_i(k))) = \inf_{ \zeta : \zeta (V_i (k)) = 
\eta (V_i (k))} p_i(a | \zeta) , $$
$$ \Delta_i^{[k]} (a |  \eta (V_i (k)))= 
r_i^{[k]} (a |\eta(V_i (k))) - r^{[k-1]}_{i} ( a | \eta ( V_i ({k-1})  )). $$ 

Then we have that
$$ p_i(a|\eta) = \sum_{j=-1}^k \Delta^{[j]}_i (a| \eta(V_i(j))) +
\left[ p_i(a|\eta) - r_i^{[k]} (a|\eta (V_i(k)))\right].
$$
By continuity of $c_i(a,\eta),$ hence of $p_i (a|\eta),$
$$ r_i^{[k]}  (a|\eta (V_i(k))) \to p_i ( a| \eta) \mbox{ as } k \to \infty .$$
Hence by monotone convergence, we conclude that 
$$ \sum_{j=-1}^\infty \Delta_i^{[j]}  (a| \eta(V_i (j))) = p_i(a| \eta )  .$$
Now, put 
$$ \lambda_i (k,\eta (V_i(k))) = \sum_a \Delta_i^{[k]} (a| \eta(V_i (k)))$$
and for any $i,k$ such that $ \lambda_i (k,\eta (V_i(k))) > 0,$ we define
$$ \tilde{p}_i^{[k]} (a | \eta(V_i(k))) = \frac{\Delta_i^{[k]} (a|
\eta(V_i (k)))}{ \lambda_i (k,\eta (V_i(k)))}.$$ For $i, k$ such that
$ \lambda_i (k,\eta (V_i(k))) = 0,$ define $\tilde{p}_i^{[k]} (a |
\eta(V_i(k)))$ in an arbitrary fixed way.

Hence
\begin{equation}\label{eq:almost}
 p_i(a|\eta) = \sum_{k=-1}^\infty \lambda_i (k,\eta(V_i(k)))
 \tilde{p}_i^{[k]} (a| \eta(V_i(k))).
\end{equation}
In (\ref{eq:almost}) the factors $\lambda_i (k,\eta(V_i(k) ))$ still
depend on $\eta (V_i(k)) .$ To obtain the decomposition as in the
theorem, we must rewrite it as follows.

For any $i,$ take $M_i$ as in (\ref{eq:m}) and the sequences $\alpha_i
(k), \lambda_i (k), k \geq -1,$ as defined in (\ref{eq:alpha}) and
(\ref{eq:lambdak}), respectively. Define the new quantities
$$\alpha_i (k,\eta(V_i(k))) = M_i \, \sum_{l \le k} \lambda_i (l, \eta( V_i(l))).$$ 

Finally put $p_i^{[-1]} (a) = \tilde{p}_i^{[-1]} (a),$ and for any $k \geq 0,$
\begin{eqnarray*}
&& p_i^{[k]} ( a| \eta(V_i({k}))) = \\
&& \sum_{-1 = l' \le l }^{k-1}  1_{\{
\alpha_i (l' - 1 ,\eta(V_i (l' -1))) < \alpha_i ({k-1}) \le \alpha_i({l'},
\eta(V_i({l' })))\}} 1_{\{
\alpha_i (l,\eta(V_i (l))) < \alpha_i ({k}) \le \alpha_i({l+1},
\eta(V_i({l+1})))\}}   \\
 && \quad \quad \quad
\left[ \frac{\alpha_i (l',\eta (V_i(l'))) - \alpha_i(k-1) }{M_i \,
\lambda_i({k})} \tilde{p}_i^{[l']} (a | \eta (V_i(l')))\right.  \\
 && \quad \quad \quad + \sum_{m = l'+1}^{l} \frac{\lambda_i (m , \eta (V_i (m))}{M_i \lambda_i (k) }
 \tilde{p}_i^{[m]} (a | \eta (V_i(m)))
 \\
 && 
\quad \quad \quad \left. + \; \frac{\alpha_i({k}) - \alpha_i (l,\eta (V_i(l)))}{M_i \,
\lambda_i({k})} \tilde{p}_i^{[l+1]} (a| \eta (V_i({l+1}))) \right] .
\end{eqnarray*} 
This concludes our proof. \hfill $\square$

\section{The black and white time-reverse sketch process} \label{sec:bw}

  The {\it black and white time-reverse sketch process} gives the
  mathematically precise description of the backward black and white
  Algorithm 1 given above. We start
  by introducing some more notation.  For each $ i \in \Z^d, $ denote
  by $\ldots T_{-2}^i <T_{-1}^i < T_{0}^i < 0 < T_1^i < T_2^i <
  \ldots$ the occurrence times of the rate $M_i$ Poisson point process
  $N^i $ on the real line. The Poisson point processes associated to
  different sites are independent. To each point $T_n^i$ associate an
  independent mark $K^i_n$ according to the probability distribution
  $(\lambda_i(k))_{k \ge -1}$. As usual, we identify the Poisson point
  processes and the counting measures through the formula
$$N^i[s,t] \,=\, \sum_{n \in \Z} \one_{\{ s \le T_n^i \le t\}}.$$
It follows from this identification that for any $t > 0$ we have
$T^i_{N^i(0,t]} \le t < T^i_{N^i(0,t]+1},$ and for any $t \le 0,$
$T^i_{-N^i(t,0]} \le t < T^i_{-N^i(t,0]+1}$.

For each $i \in \Z^d$ and $t \in \R$ we define the time-reverse point
process starting at time $t,$ associated to site $i,$
\begin{eqnarray}
  \label{eq:tildet}
  \tilde{T}^{(i,t)}_n &=& t \,-\, T^i_{N^i(0,t]-n+1}, \quad t \ge 0,\nonumber \\
  \tilde{T}^{(i,t)}_n & = &  t \,-\, T^i_{-N^i(t,0]-n+1}, \quad t < 0 .
 \end{eqnarray}
 We also define the associated marks
\begin{eqnarray}
  \label{eq:tildet1}
  \tilde{K}^{(i,t)}_n &=& K^i_{N^i(0,t]-n+1}, \quad t \ge 0,\nonumber \\
  \tilde{K}^{(i,t)}_n & = &  K^i_{-N^i(t,0]-n+1}, \quad t < 0.
 \end{eqnarray}

 For each site $i \in \Z^d$, $k \ge -1$, the reversed $k$-marked
 Poisson point process returning from time $t$ is defined as
 \begin{equation}
   \label{eq:tilden}
   \tilde{N}^{(i,t,k)}[s,u] \,=\, \sum_{n} \one_{\{s \le \tilde{T}^{(i,t)}_{n} \le u\}} \one_{\{\tilde{K}^{(i,t)}_n = k\}}.
 \end{equation}

 To define the black and white time-reverse sketch process we need to
 introduce a family of transformations $\{\pi^{(i,k)}, i \in \Z^d, k
 \ge 0\}$ on the set of finite subsets of $\Z^d,$ $ {\cal F}(\Z^d),$
 defined as follows. For any unitary set $\{j\}$,
 \begin{equation}
   \label{eq:pij}
   \pi^{(i,k)}(\{j\}) \,=\, \left\{ \begin{array}{ll}
                                   V_i (k), & \mbox{ if } j=i \\
                                   \{j\}, & \mbox{ otherwise}
                                    \end{array} \right\} .
 \end{equation}
 Notice that for $k=-1,$ $ \pi^{(i,k)}(\{i\}) = \emptyset.$ For any
 set finite set $F \subset \Z^d$, we define similarly
\begin{equation}
  \label{eq:pif}
  \pi^{(i,k)}(F) \,=\, \cup_{j \in F} \pi^{(i,k)}(\{j\}) . 
\end{equation}

The black and white time-reverse sketch process starting at site $i$
at time $t$ will be denoted by $(C_s^{(i,t)})_{s \geq 0}.$
$C_s^{(i,t)}$ is the set of sites at time $s$ whose colors affect the
color of site $i$ at time $t.$ The evolution of this process is
defined through the following equation: $C_0^{(i,t)} := \{i\},$ and
\begin{equation}
  \label{eq:ct}
 f( C_s^{(i,t)}) \,=\, f(C_0^{(i,t)}) \,+\, \sum_{k \ge -1} \sum_{j \in \Z^d} \int_0^s [f(\pi^{(j,k)} (C_{u-}^{(i,t)})) - f(C_{u-}^{(i,t)})]\, \tilde{N}^{(j,t,k)}(du), 
\end{equation}
where $f: {\cal F}(\Z^d) \rightarrow \R$ is any bounded cylindrical
function. This family of equations characterizes completely the time
evolution $\{C_s^{(i,t)}, s \ge 0\}$. For any finite set $F \subset
\Z^d$ define
$$C_s^{(F,t)} \,=\, \cup_{i \in F} C_s^{(i,t)}.$$

The following proposition summarizes the properties of the family of
processes defined above.

\begin{prop}
  For any finite set $F \subset \Z^d$, $C_s^{(F,t)}$ is a Markov jump
  process having as infinitesimal generator
\begin{equation}
  \label{eq:generatord}
  L f(C) \,=\, \sum_{i \in C} \sum_{k \ge 0} \lambda_i (k) [f(C  \cup V_i(k)) - f(C)] + \lambda_i (-1) [f(C \setminus \{i\}) - f(C)] , 
\end{equation}
where $f$ is any bounded function.
\end{prop}

{\bf Proof:} The proof follows in a standard way from the construction
\reff{eq:ct}.

\section{Proof of Theorem \ref{theo:1}}
The existence issue addressed by Theorem \ref{theo:1} can be
reformulated in terms of the black and white time-reverse sketch
process described above. The process is well defined if for each site
$i$ and each time $t$, the time-reverse procedure $C^{(i,t)} $
described above is a non-explosive Markov jump process. This means
that for each time $t,$ the number of operations needed to determine
the value of $\xi^{\eta}_t(i)$ is finite almost surely. Note that by
equation (\ref{eq:ct}), the jumps of $C^{(i, t)} $ occur at total rate
  $$ \sum_{j \in C_s^{(i,t)}}  M_j \, \sum_{ k \geq - 1 } \lambda_j (k) \le (\sup_j M_j ) \, | C_s^{(i,t)} |,$$
  where $| \cdot | $ denotes the cardinal of a set. Hence it suffices to show that the cardinal of 
  $C^{(i, t) }$ remains finite. 

  More precisely, fix some $N \in \N.$ Let $L_s = | C_s^{(i,t)}| $ and
$$ T_N = \inf \{ t : L_t \geq N \} .$$ Then by (\ref{eq:ct}),
\begin{eqnarray}\label{eq:ub}
  L_{s \wedge T_N} & \le&  1 + \sum_{k \geq 1} \sum_{j \in \Z^d} \int_0^{s \wedge T_N} [|V_j(k)| - 1] 1_{\{ j \in C^{(i,t)}_{u-}  \}}  \, \tilde{N}^{(j,t,k)}(du)\nonumber\\
  &&  -\sum_{j \in \Z^d} \int_0^{s \wedge T_N}  1_{\{ j \in C^{(i,t)}_{u-}  \}}  \, \tilde{N}^{(j,t,0 )}(du)  .
 \end{eqnarray}
 Passing to expectation and using that by condition
 (\ref{eq:condition1}),
$$ m = \sup_i \sum_{k \geq 1} M_i \, \lambda_i (k)  |V_i(k)|  < + \infty,$$
this yields
\begin{eqnarray}\label{eq:upperbound2}
E (L_{s \wedge T_N})& \le& 1 +  \sum_{j \in \Z^d}\, M_j  \left( (\sum_{k \geq 1}\, \lambda_j (k) [ |V_j(k)| - 1] ) -  \lambda_j (-1)\right) \nonumber \\
&& \quad \quad \quad \quad \quad \quad  \quad \quad \times  E \int_0^{s\wedge T_N} 1_{\{ j \in C^{(i,t)}_{u-}  \}}  du  \nonumber \\
&\le& 1 +  m \, E \int_0^{s\wedge T_N}  L_u du .
\end{eqnarray}
Letting $N \to \infty ,$ we thus get that
$$ E(L_s) \le 1 + m \int_0^s E(L_u) du,$$
and Gronwall's lemma yields 
\begin{equation}
\label{eq:gron} 
 E(L_s) \le e^{ms} .
\end{equation}

This implies that the number of sites that have to be determined in
order to know the value of site $i$ at time $t$ is finite almost
surely.  This means that the process  $C_s^{(F, t)}
$ admits only a finite number of
jumps on any finite time interval. Hence, we have necessarily
$N_{STOP} < + \infty $ almost surely which means that the algorithm
stops almost surely after a finite time. This concludes the proof of Theorem \ref{theo:1}.

\section{Proof of Theorem \ref{theo:nstop}}
We show that under condition (\ref{eq:condition2}),  
Algorithm 3 stops after a finite
time almost surely. Write $L^i_s$ for the cardinal of
$C_s^{(i, t)} .$ Using once more the upper-bound (\ref{eq:upperbound2})
and the fact that under condition (\ref{eq:condition2}),
$$ M_j \left( (\sum_{k \geq 1}\, \lambda_j (k) [ |V_j(k)| - 1] ) -  \lambda_j (-1)\right)  \le - \varepsilon < 0 ,$$
Gronwall's lemma yields that
$$ E( L^i_s) \le e^{- \varepsilon s} .$$ 
Hence, since $ |C_s^{(F, t)} | \le \sum_{i \in F } |C_s^{(i, t)}| =
\sum_{i \in F } L_s^i,$
$$ E ( |C_s^{(F, t)} |)  \le |F|  e^{ - \varepsilon s .} $$
This implies that $ \inf \{ s: C_s^{(F, t)} = \emptyset \} $ is finite
almost surely. Due to Theorem \ref{theo:1}, the process $C_s^{(F, t)}
$ is non-explosive, which means that it admits only a finite number of
jumps on any finite time interval. Hence, we have necessarily
$N_{STOP} < + \infty $ almost surely which means that the algorithm
stops almost surely after a finite time.

In order to show that the measure $\mu $ that we have simulated in this
way is necessarily the unique invariant probability measure of the
process, we prove the following lemma.

\begin{lem}\label{lemma:1}
Fix a time $t > 0,$ some finite set of sites $F \subset \Z^d $ and two
initial configurations $ \eta $ and $\zeta \in A^{\Z^d} .$ Then there exists a
coupling of the two processes $(\xi^\eta_s)_s$ and $(\xi^\zeta_s)_s $ such that
$$ 
P( \xi_t^\eta (F) \neq \xi^\zeta_t (F) ) \le |F| e^{- \varepsilon t } .
$$  
\end{lem} 

From this lemma, it follows immediately that $\mu$ is the unique
invariant measure of the process and that the convergence towards the
invariant measure takes place exponentially fast.

{\bf Proof of Lemma \ref{lemma:1}}.
We use a slight modification of Algorithm 1 and 2 in order to 
construct $\xi_t^\eta $ and $\xi_t^\zeta . $ The modification is defined
as follows. Replace step 8 of Algorithm 1 by 
$$ \mbox{ \bf if } K = -1 ,\mbox{ \bf then } $$
$$ C \leftarrow C \setminus \{ I \}$$
$$ \mbox{ \bf else } $$
$$ C \leftarrow C \cup V_I (K) $$ We use the same realizations of $T,
I $ and $ K $ for the construction of $\xi_t^\eta $ and $\xi_t^\zeta
$.  Write $L_s$ for the cardinal of $C_s^{(F,t)}.$ Clearly, both
realizations of $\xi_t^\eta $ and $\xi_t^\zeta $ do not depend on the
initial configuration $\eta,$ $\zeta$ respectively if and only if the
output $C$ of Algorithm 1 is void. Thus,
\begin{eqnarray*}
P ( \xi_t^\eta   ( F) \neq \xi_t^\zeta (F) )
 & \le &  P( T_{STOP}  \geq  t )  \\
 & =& \P ( L_t \geq 1 )\\
 & \le & E (L_t) \le |F| e^{- \varepsilon t }   .
\end{eqnarray*}
This concludes the proof of lemma \ref{lemma:1}. 

\section{Proof of Theorem \ref{theo:4}} 

The proof of Theorem \ref{theo:4} goes according to the following
lines. 

We start at time $t$ with the sites in $F$ and go back
into the past until time 0 following the backward black and white
sketch without deaths described in Algorithm 1. The set $C$ of points
reached by this procedure at time 0 are the only ones which coloring
affects the evolution during the interval of time $[0,t]$ of the sites
belonging to $F$.

We need to know that $C$ is a finite set. This follows from a slight
modification of the proof of Theorem \ref{theo:1}. Notice that in the
construction by Algorithm 1, even if $K = -1 ,$ the corresponding site
is not removed from the set $C. $ This implies that in the upper bound
(\ref{eq:ub}) the negative term on the right hand side disappears.
This modification does not affect the upper bound
(\ref{eq:upperbound2}) which remains true.

Using Theorem \ref{theo:nstop} we assign colors to the sites
belonging to $C$ using the invariant distribution at the origin of the
multicolor system. Then we apply Algorithm 2 to describe the time evolution
of the coloring of the sites in $F$. This evolution depends on the
colors of the sites in $C$ at time zero as well as on the successive
choices of sites and ranges made during the backward steps
starting at time $t$. This concludes the proof.

\section{Proof of Theorem \ref{theo:5} }

The construction of the finitary coding can be better understood if we
do it using two intermediate steps based on families of
infinite-valued random variables.

Using a slightly abusive terminology, let us introduce the following
definitions of a finitary coding from families of piles of
i.i.d. random variables.

\begin{defin}\label{def:finitarycoding}
  We say that there exists a {\it finitary coding} from a family of
  i.i.d. uniform random variables $(U_n (i) ,i
  \in \Z^d, n \in \N )$ to the configuration $(\xi_0 (i), i \in \Z^d)
  $ sampled with respect to the invariant probability measure $\mu$,
  if there exists a deterministic function $f:[0,1]^{\Z^d \times \N}
  \rightarrow A^{\Z^d}$ such that, almost surely, the following holds
\begin{itemize}
\item $f$ commutes with the shift operator in $\Z^d$;
\item $\xi_0 = f( (U_n (j) ), j \in \Z^d, n \in \N  ) $; and
\item for each site $i \in \Z^d$, there exists a a finite subset
  $\bar{F}_i $ of $\Z^d$ and $\bar{n}_i \ge 1$ such that if
$$ U'_n (j) = U_n (j) \mbox{ for all } j \in \bar{F}_i,  n \le
\bar{n}_i$$
then
$$  f ( (U_n (j) ), j \in \Z^d, n \in \N)(i) =  f ( (U'_n (j) ), j \in
\Z^d, n \in \N)(i).$$
\end{itemize}
\end{defin}

\begin{defin}\label{def:finitarycoding2}
  We say that there exists a {\it finitary coding} from a family of
  i.i.d. fair Bernoulli random variables $(Y_{n,r} (i) ,i \in \Z^d, (n,r)
  \in \N^2 )$ to the configuration $(\xi_0 (i), i \in \Z^d) $ sampled
  with respect to the invariant probability measure $\mu$, if there
  exists a deterministic function $f:\{0,1\}^{\Z^d \times \N^2}
  \rightarrow A^{\Z^d}$ such that, almost surely, the following holds
\begin{itemize}
\item $f$ commutes with the shift operator in $\Z^d$;
\item $\xi_0 = f( (Y_{n,r} (j) ), j \in \Z^d, (n,r) \in \N^2  ) $; and
\item for each site $i \in \Z^d$, there exists a a finite subset
  $\bar{F}_i $ of $\Z^d$ and $\bar{n}_i \ge 1$ such that if
$$ Y'_{n,r} (j) = Y_{n,r} (j) \mbox{ for all } j \in \bar{F}_i, 1 \le n \le
\bar{n}_i, 1 \le r \le \bar{n}_i$$
then
$$  f ( (Y_{n,r} (j) ), j \in \Z^d, (n,r) \in \N^2)(i) =  f (
(Y'_{n,r} (j) ), j \in \Z^d, (n,r) \in \N^2)(i).$$
\end{itemize}
\end{defin}

We will first prove the existence of a finitary coding in this
extended definition from the sequence $(U_n(i), i \in \Z^d , n \in \N
)$ to the configuration $(\xi_0 (i), i \in \Z^d) $ .  

\begin{prop} \label{prop:uni-fini} Under the conditions of Theorem
  \ref{theo:nstop}, there exists a finitary coding from a family of
  i.i.d. uniform random variables $(U_n (i) ,i \in \Z^d, n \in \N )$
  to the configuration $(\xi_0 (i), i \in \Z^d)$.
\end{prop}

\proof Our goal is to choose the color of site $i$ at time 0 using
Algorithms 3 and 4.  For notational convenience, we will represent the
sequence $U_n(j) $ as
$$U_n^v (j), \, v \in {\cal V} \, ,$$
where ${\cal V} = \{I, K, W\}$.

The first sequence of uniform variables $U_n^I (j) $ is used for the
choice of successive points $I$ at Step 5 of Algorithm 3.  The second
sequence $U_n^K(j)$ will be used to construct the sequence of ranges
$K$ of Step 6 of Algorithm 4. Finally, the third sequence $U_n^W(j)$ will be
used to construct the corresponding colors $W$ at Steps 7 and 9 of the forward
procedure described in Algorithm 4.

The fact that $N_{STOP } $ is finite almost surely implies that the
backward Algorithm 1 must run only a finite number of steps for
any fixed $i$. Thus only a finite set of sites $\bar{F}_i $ is
involved in this procedure. This also implies that the number of
uniform random variables we must use is finite, thus $\bar{n}_i$ is
finite.  The definition of the function $f$ is explicitly given by
Algorithms 3 and 4. 

This concludes the proof of Proposition \ref{prop:uni-fini}.

Proposition \ref{prop:uni-fini} can be rewritten using piles of piles
of Bernoulli random variables instead of piles of uniform random variables.

\begin{prop} \label{prop:bern-fini} Under the conditions of Theorem
  \ref{theo:nstop}, there exists a finitary coding from a family of
  i.i.d. fair Bernoulli random variables $(Y_{n,r} (i) ,i \in \Z^d,
  (n,r) \in \N^2 )$ to the configuration $(\xi_0 (i), i \in \Z^d)$.
\end{prop}

\proof  As before,  for notional convenience, we will represent the
sequence $Y_{n,r}(j) $ as $$Y_{n,r}^v (j), \, v \in {\cal V} .$$

All we need to prove is that the successive uniform random variables
$\{U^v_n(j), j \in \bar{F}_i, 1 \le n \le \bar{n}_i\}$ used in Proposition 
\ref{prop:uni-fini} can be generated
using a finite number of Bernoulli random variables from the pile
$Y_{n,r}^v (j), r \in \N$.

In all the steps, the uniform random variables were used to generate 
random variables taking values in a countable set. Let us identify this
countable set with $\N$. In the successive steps of Algorithms 3 and
4, this selection could generate either $K$, $I$ or $W$. The selection
is made by defining in each case a suitable partition of $[0, 1 ] =
\cup_{l=1}^{\infty} [\theta(l), \theta(l+1))$ and then choosing the
index $l$ whenever $U^v_n(j) \in [\theta(l), \theta(l+1))$. It is easy to see
that $U_n^v(j)$ has the same law as $\sum_{r=1}^{\infty} 2^{-r}
Y_{n,r}^v (j)$. 

We borrow from Knuth and Yao (1976) the following algorithm to
generate the discrete random variables $I$, $K$ and $W$ which appear
in the perfect sampling procedure (see also Harvey {\it et al.}, 2005).   
For any $m \geq 1,$ we define
$$ S_m (j,n,v) = \sum_{r=1}^m 2^{-r} Y_{n,r}^v (j) .$$
Now, put 
$$ J(S_m(j,n,v)) = \sup \{ k \ge 1 : \theta(k) \le S_m(j,n,v) \} ,$$
and finally define 
\begin{equation}
\label{eq:N(v,n,j)}
 N_n^v(j) = \inf \{ m \ge 1 : J(S_m(j,n,v)) = J(S_{m'}(j,n,v) )\,
 \forall m' \geq m \} .
\end{equation}

Notice that $ N_n^v(j)$ is a finite stopping time with respect to the
$\sigma$-algebra generated by $\{Y_{n,r}^v(j), r \ge 1 \}$.
Therefore, the total number of piles $N(j)$ is equal to $ N(j)= \sum_{v \in {\cal V}}
\sum_{n=1}^{\bar{n}_j}N_n^v(j), $ where $\bar{n}_j$ is defined in
Proposition \ref{prop:uni-fini}, used at site $j$ is finite and the
event $[N(j) = \ell]$ is measurable with respect to the
$\sigma$-algebra generated by $\{Y_{n,r}^v(j), v \in {\cal V}, 1 \le n
\le \ell, 1 \le r \le \ell\}$. Since the set of sites used is $\bar{F}_i$
(the same one as Proposition \ref{prop:uni-fini}), the proof is
complete.

  We are finally ready to prove Theorem \ref{theo:5}.  The difficulty
  is to show that the construction achieved in Proposition
  \ref{prop:bern-fini} using the random sized piles $(Y_{n,r}^v (i) ,
  v\in {\cal V} , i \in \Z^d, n \le N(i), r \le \ N(i) )$ can actually be
  done using finite piles of fair Bernoulli random
  variables. Specifically, for $i \in \Z^d ,$ let us call 
  $$Z(i) = (Y_{n,r}^v (i) , v\in {\cal V},  n \in \{0,\ldots, M\}, r \in
  \{0,\ldots, M\} ),$$
 where $M$ is a suitable fixed positive integer.
 
 The proof that there exists a finitary coding from the family of
 finite-valued i.i.d. random variables $\{Z(i), i \in \Z^d\}$ to
 $\{\xi_0(i), i \in \Z^d\}$ follows from the construction in Van den
 Berg and Steif (1999) if we can take
$$ M > \sup_{j \in \Z^d} \E[N(j)]. $$

It follows from the definition of $N_n^v(j)$ given by
(\ref{eq:N(v,n,j)}) that
\begin{equation}
 \label{eq:tail1}
\P[ N_n^v(j) > k]  \le  \P \left( \cup_{i=1}^{m_k}  \left[ \theta(i)
    - \frac{1}{2^k} < S_k (j,n,v) \le \theta(i) \right] \right) + \P \left(  1
    - \frac{1}{2^{k}} < S_k(j,n,v) \le  1  \right)
\end{equation}
where
$$m_k = \sup\{ i \ge 1; \theta(i) < 1 - \frac{1}{2^k} \}.$$
 
In the above formula,  the partition of $[0, 1 ] =
\cup_{i=1}^{\infty} [\theta(i), \theta(i+1))$ was used to
simulate the countable-valued random variable at stake at that level
(either $I$, $K$ or $W$). 

Since $ S_k (j,n,v)$ converges in law to a uniform random variable as
$k \rightarrow \infty$, the right hand side of (\ref{eq:tail1}) is
bounded above by $$\frac{m_k+1}{2^{k-1}}.$$

For piles choosing colors, the result is obvious since the set of
possible colors is finite and therefore all the corresponding $m_k =
|A|$ for all $k \ge |A|$. 

For piles choosing sites in the backward black and white sketch, 
the result follows from  inequality (\ref{eq:gron}) as in the conclusion of the proof of Proposition \ref{theo:1}. 

Finally, for piles choosing ranges, the result follows from the following
two lemmas and Wald's inequality observing that  $ \bar{n}_j \le N_{STOP}$.

\begin{lem} $\sup_j \E[N_n^K(j)] < \infty.$ \end{lem}

\proof It follows from Knuth and Yao (1976) that
\begin{equation}
  \label{eq:entropy}
  \E[N_n^K(j)] \le H(\{ \lambda_j(k) \}_{k \ge -1}) + 2,
\end{equation}
where $H(\{ \lambda_j(k)\}_{k \ge -1 })$ is the entropy of the discrete
distribution $\{ \lambda_j(k)\}_{k \ge -1 }$ of the random variable $K$.
On the other hand, the condition $\sum_k |V_k(j)| \lambda_k(j)  < 1$ in
Theorem \ref{theo:nstop} implies that $\sum_k k \lambda_j(k) = m_j < 1$
for all $j \in \Z^d$. We want to compare $\{ \lambda_j(k)\}_{k \ge -1 }$
to a geometric distribution. For that sake, we introduce a distribution $\nu_j (k) 
$ on $\{ 1, 2, \ldots \} $ by $ \nu_j (k) = \lambda_j (k-2) .$ Then 
$$ \sum_{k \geq 1 } k \nu_j (k) = m_j + 2 - \lambda_j (-1) =: \tilde m_j.$$
By a direct comparison with the geometric
distribution of mean $\tilde m_j$ we have that 
\begin{equation}
  \label{eq:geom}
  H(\{ \lambda_j(k)\}_{k \ge -1}) \le - \log(p_j) - \log (1-p_j) (\tilde m_j-1) < \infty,
\end{equation}
where $p_j = 1/(\tilde m_j)  $. 

\begin{lem} $\E[N_{STOP}] < \infty$. \end{lem}

\proof Without loss of generality we can consider $F = \{0\}$ to start
Algorithm 3. Define
$$ L_n := | C_n|,$$
the cardinal of the set $C_n$ after $n$ steps of Algorithm 3. Let
$(K^{i}_n)_{n \geq 0, i \in \Z^d} $ be the i.i.d. marks defined in
Section \ref{sec:bw}, taking values in $\{ -1 , 0 , 1, 2, \ldots
\} $ such that 
$$ P( K^i_n = k ) = \lambda_i (k) .$$
Define $X^i_n = |V_0(K_n^i) | - 1$. Note that by condition
(\ref{eq:condition2}), 
$$\sup_{i \in \Z^d } E(X^i _1) \le (\bla -1) <
0,$$
where $\bla =  \sup_{i \in \Z^d}   \sum_{k \ge 0} \, |V_i (k)|
\lambda_i (k) $. 

Consider the sequence $I_n $ which gives the site of the particle
chosen at the $n$th step of Algorithm 3. Put
$$S_n \,=\, \sum_{k= 0}^n X_k^{I_k}.$$
Note that by construction, $S_n + n (1 - \bla) $ is a super-martingale. 
Then a very rough upper bound is 
$$ L_n \le 1 + S_n \mbox{ as long as } n \le V_{STOP} ,$$
where $V_{STOP} $ is  defined as 
$$ V_{STOP} = \min \{ k : S_k = -1 \} .$$
By construction
$$  N_{STOP} \le V_{STOP}.  $$
Fix a truncation level $N.$ Then by the stopping rule for super-martingales, we have that
$$ E( S_{V_{STOP} \wedge N } ) + (1 - \bla) E( V_{STOP} \wedge N ) \le 0  .$$
But notice that 
$$  E( S_{V_{STOP} \wedge N }) = - 1 \cdot  P( V_{STOP} \le N ) + E (S_N ; V_{STOP } > N ) .$$   
On $   V_{STOP } > N , $ $S_N \geq 0,$ hence we have that $   E( S_{V_{STOP} \wedge N }) \geq - P( V_{STOP} \le N) .$ 
We conclude that 
\begin{eqnarray*}
E( V_{STOP} \wedge N ) & \le & \frac{1}{(1 - \bla) }  P( V_{STOP} \le N ) .
\end{eqnarray*}
Now, letting $N \to \infty ,$ we get 
$$ E( V_{STOP}) \le \frac{1}{(1-\bla)} ,$$ 
and therefore
$$ E( N_{STOP}) \le \frac{1}{(1-\bla)}. $$ 

\vspace{.5cm}

Finally for piles choosing sites in the backward black and white
sketch, the result follows from  inequality (\ref{eq:gron}) as in
the conclusion of the proof of Theorem \ref{theo:1}. 

This concludes the proof.

\section{Proof of Corollary \ref{cor:HS}}

 The strategy of the proof is the following. We will consider a
multicolor system having the law of the Markov random field as
invariant measure. Then the corollary follows from the second
assertion of Theorem \ref{theo:5}.  

A standard way to obtain a system having the law $\mu$ as invariant
measure is to ask for reversibility. Usually, in the statistics
literature such dynamic is known as Gibbs sampler. In the statistical
physics literature, where it first appeared, it is known as Glauber
dynamics. We will use a particular case of the Glauber dynamics called
the heat bath algorithm. The idea is that at each site there is a
Poisson clock which rings independently of all other sites. Each
time its clock rings the color of the site is updated according to
the specification of the Markov random field ${\cal Q}$.

We are only considering the Markov spatially homogeneous
case. This means that the rate $c_i ( a , T_i \xi ) = c_0 ( a, \xi )$
where $(T_i \xi ) (j) = \xi ( j-i ) .$ $c_0(a, \xi) $ only depends on
$\xi (\partial 0  ) $. Moreover, in this case $ M = M_i = 1$.

With the heat bath algorithm, the rates are defined as 
$$c_0(a, \xi) = {\cal Q}(X(0)=a|
X(\partial 0) = \xi(\partial 0)) $$ where ${\cal Q}$ is the {\it
  specification} of the random field $X$. 

Since we are considering the homogeneous case, we will drop the
subscript from the notation. Therefore, we have
\begin{eqnarray*}
  \alpha(-1) &=&  \sum_{a \in A} \min \left( \inf_{\zeta \in A^{\Z^d},
      \zeta(i ) \neq a} c_0(a, \zeta) 
    , \, 1 - \sup_{ \zeta \in A^{\Z^d}, \zeta(i )= a} \sum_{b \neq a}
    c_0 (b, \zeta) \right)  \\
 &=& \sum_{a \in A} \min \left( \inf_{\zeta \in A^{\Z^d},
      \zeta(i ) \neq a} c_0(a, \zeta) 
    , \, \inf_{ \zeta \in A^{\Z^d}, \zeta(i )= a} c_0 (a, \zeta)
  \right) \\
& =& \sum_{a \in A} \min_{\zeta(\partial 0) \in A^{\partial 0}} {\cal Q}(X(0)=a|
X(\partial 0) = \zeta(\partial 0)).  \\
\end{eqnarray*}
Also,
\begin{eqnarray*}
  \alpha (0) &=& \min_{w \in A } \left( \sum_{a \in A, a \neq w}
  \inf_{\zeta: \zeta(0) = w} c_0( a,\zeta) + 1 - \sup_{\zeta:
    \zeta(0) = w} \sum_{ b \neq w} c_0 (b, \zeta)\right)  \\
&=& \min_{w \in A } \sum_{a \in A}  \inf_{\zeta: \zeta(0) = w} c_0( a,\zeta) \\
&=& \min_{w \in A }  \sum_{a \in A}  \inf_{\zeta;\zeta(0) = w} {\cal Q}(X(0)=a|
X(\partial 0) = \zeta(\partial 0)) \\
&=&  \sum_{a \in A} \min_{\zeta(\partial 0) \in A^{\partial 0}} {\cal Q}(X(0)=a|
X(\partial 0) = \zeta(\partial 0)).
\end{eqnarray*}
The last equality follows from the fact that $ \inf_{\zeta;\zeta(0) =
  w} {\cal Q}(X(0)=a| X(\partial 0) = \zeta(\partial 0))$ only 
depends on the value of the random field at $\partial 0$. 

Finally, for all $k \ge 1$
\begin{eqnarray*}
  \alpha (k) &=& \min_{w \in A^{V_0 (k) }} \left( \sum_{a \in A, a \neq w(0)}
    \inf_{\zeta: \zeta(V_0 (k)) = w} c_0( a,\zeta) + 1 - \sup_{\zeta:
      \zeta(V_0 (k)) = w} \sum_{ b \neq w(0)} c_0 (b, \zeta)\right) \\
  &=&  \min_{w \in A^{V_0 (k) }} \left( \sum_{a \in A}
    \inf_{\zeta: \zeta(V_0 (k)) = w} c_0( a,\zeta) \right) \\
  &=&  \min_{w \in A^{\partial 0 }} \left( \sum_{a \in A}
    {\cal Q}(X(0)=a|
    X(\partial 0) = w(\partial 0) \right) = 1.
\end{eqnarray*} 

Observe that $\alpha(0) = \alpha(-1)$. Therefore, condition
(\ref{eq:condition1}) reduces to
$$\alpha(-1) >  \frac{2d}{2d +1 }.$$

\section*{Acknowledgments}

We thank Pablo Ferrari, Alexsandro Gallo, Yoshiharu Kohayakawa, Servet
Martinez, Enza Orlandi and Ron Peled for many comments and
bibliographic suggestions. We also thank the anonymous Associated
Editor that pointed out an incomplete definition in an earlier version of
this manuscript.

This work is part of PRONEX/FAPESP's project \emph{Stochastic
behavior, critical phenomena and rhythmic pattern identification in
natural languages} (grant number 03/09930-9), CNRS-FAPESP project
\emph{Probabilistic phonology of rhythm} and CNPq's projects
\emph{Stochastic modeling of speech} (grant number 475177/2004-5) and
\emph{Rhythmic patterns, prosodic domains and probabilistic modeling
in Portuguese Corpora} (grant number 485999/2007-2). AG and NLG are
partially supported by a CNPq fellowship (grants 308656/2005-9 and
301530/2007-6, respectively).

\def\refname{References}

\vskip30pt

Antonio Galves

Instituto de Matem\'atica e Estat\'{\i}stica

Universidade de S\~ao Paulo

PO Box 66281

05315-970 S\~ao Paulo, Brasil

e-mail: {\tt galves@ime.usp.br}

\bigskip

Nancy L. Garcia

Instituto de Matem\'atica, Estat\'{\i}stica e Computa\c c\~ao
Cient\'\i fica 

Universidade Estadual de Campinas

PO Box 6065

13083-859 Campinas, Brasil

e-mail: {\tt nancy@ime.unicamp.br}

\bigskip

Eva L\"ocherbach

Universit\'e Paris-Est

LAMA -- UMR CNRS 8050

61, Avenue du G\'en\'eral de Gaulle 

94000 Cr\'eteil, France

e-mail: {\tt locherbach@univ-paris12.fr}

\end{document}